\documentclass{amsart}
\usepackage{amssymb,algorithm,algorithmic,url,hyperref}

\newtheorem{theorem}{Theorem}[section]
\newtheorem{lemma}[theorem]{Lemma}
\newtheorem{corollary}[theorem]{Corollary}
\theoremstyle{definition}
\newtheorem{definition}[theorem]{Definition}

\newcommand{\la}{\langle}
\newcommand{\ra}{\rangle}

\newcommand{\ZZ}{\mathbb{Z}}  
\newcommand{\Aut}{{\rm {Aut}}}
\newcommand{\Cos}{{\rm {Cos}}}

\newcommand{\vW}{\vec{{\rm W}}}

\newcommand{\Pl}{{\rm Pl}}
\newcommand{\alt}{\mathcal{A}}
\newcommand{\V}{{\rm V}}
\newcommand{\A}{{{\rm A}}}
\newcommand{\E}{{{\rm E}}}

\newcommand{\val}{\mathop{\rm val}}
\newcommand{\hal}{\frac{1}{2}}
\newcommand{\C}{{\mathcal{C}}}
\newcommand{\D}{{\mathcal{D}}}
\newcommand{\N}{{\mathcal{N}}}
\newcommand{\cA}{{\mathcal{A}}}
\newcommand{\tG}{{\tilde{G}}}
\newcommand{\am}{A}
\newcommand{\half}{$\frac{1}{2}$}
\newcommand{\oppG}{{\Gamma^{\rm{opp}}}}

\begin{document}

\title[Census of $4$-valent half-arc-transitive and arc-transitive graphs]{A census of $4$-valent half-arc-transitive graphs and arc-transitive digraphs of valence two}
{} 

\author{Primo\v z Poto\v cnik}            
\address{Faculty of Mathematics and Physics, University of Ljubljana, \\ Jadranska 19, 1111 Ljubljana, Slovenia\\ and \\ IAM, University of Primorska,\\ Muzejski trg 2, 6000 Koper, Slovenia} 
\email{primoz.potocnik@fmf.uni-lj.si}

\author{Pablo Spiga}            
\address{University of Milano-Bicocca, Dipartimento di Matematica Pura e Applicata, \\
 Via Cozzi 53, 20126 Milano, Italy}
\email{pablo.spiga@unimib.it}

\author{Gabriel Verret}            
\address{Centre for Mathematics of Symmetry and Computation, The University of Western Australia, 35 Stirling Highway, Crawley, WA 6009, Australia\\ and\\ FAMNIT, University of Primorska, Glagolja\v{s}ka 8, SI-6000 Koper, Slovenia} 
\email{gabriel.verret@uwa.edu.au}

\thanks{The first author is supported by Slovenian Research Agency, projects L1--4292 and P1--0222. The third author is supported by UWA as part of the Australian Research Council grant DE130101001.\\ Corresponding author. Supported by Slovenian Research Agency, projects L1--4292 and P1--0222.}

\keywords{graph, digraph, edge-transitive, vertex-transitive, arc-transitive, half-arc-transitive}               
\subjclass[2010]{05E18, 20B25}                       

\begin{abstract}
A complete list of all connected arc-transitive asymmetric digraphs of in-valence and out-valence $2$ on up to $1000$ vertices
is presented. As a byproduct, a complete list of all connected $4$-valent graphs admitting a $\frac{1}{2}$-arc-transitive group of automorphisms on up to $1000$ vertices
is obtained. Several graph-theoretical properties of the elements of our census are calculated and discussed.
\medskip
\begin{center}
{\em Dedicated to Dragan Maru\v si\v c on the occasion of his 60th birthday}
\\ ${}$ \\ 
\end{center}
\end{abstract}
\maketitle

\section{Introduction}
\label{sec:intro}

Recall that a graph $\Gamma$ is called {\em $\hal$-arc-transitive} provided that its automorphism group $\Aut(\Gamma)$ acts transitively
on its edge-set $\E(\Gamma)$ and on its vertex-set $\V(\Gamma)$ but intransitively on its arc-set $\A(\Gamma)$. More generally, 
if $G$ is a subgroup of $\Aut(\Gamma)$ such that $G$ acts transitively on $\E(\Gamma)$ and $\V(\Gamma)$ but intransitively on $\A(\Gamma)$,
then $G$ is said to {\em act $\hal$-arc-transitively} on $\Gamma$ and we say that $\Gamma$ is {\em $(G,\hal)$-arc-transitive}. To shorten notation,
we shall say that a  $\hal$-arc-transitive graph is a \emph{HAT} and that a  graph admitting
a $\hal$-arc-transitive group of automorphisms is a \emph{GHAT}. Clearly, any HAT is also a GHAT. Conversely, a GHAT 
is either a HAT or arc-transitive.

The history of GHATs goes back to Tutte who, in his 1966 paper \cite[7.35, p.59]{tutte}, proved that every GHAT is of even valence and
asked whether HATs exist at all. The first 
examples of HATs were discovered a few years later by Bouwer \cite{Bou}. After a short break,
interest in GHATs picked up again in the 90s, largely due to a series of influential papers of  Maru\v{s}i\v{c}
 concerning the GHATs of valence $4$ (see \cite{AlsMarNow,Mar98,MarPra,MarXu}, to list a few).
For a nice survey of the topic, we refer the reader to \cite{MarSurvey}, and for an overview of some more recent results, see \cite{KutMarSpaWanXu,MarSpa}.

To shorten notation further, we shall say that a connected GHAT (HAT, respectively) of valence $4$ is a $4$-GHAT ($4$-HAT, respectively). 
The main result of this paper is a compilation of a complete list of all $4$-GHATs with at most $1000$ vertices. This result was obtained indirectly using
an intimate relation between $4$-GHATs and connected arc-transitive asymmetric digraphs of in- and out-valence $2$ (we shall call
such digraphs $2$-ATDs for short) -- see Section~\ref{sec:ATvsGHAT} for details on this relationship.
These results can be succinctly summarised as follows:

\begin{theorem}
There are precisely 26457 pairwise non-isomorphic $2$-ATDs on at most $1000$ vertices, and 
precisely 11941 $4$-GHATs on at most $1000$ vertices, of which 8695 are arc-transitive and 3246 are $\hal$-arc-transitive.
\end{theorem}

The actual lists of (di)graphs, together with a spreadsheet (in a ``comma separated values'' format) with some graph theoretical invariants,
is available at \cite{online}.

The rest of this section is devoted to some interesting facts gleaned from these lists.
All the relevant definitions that are omitted here can be found in Section~\ref{sec:not}. In Section~\ref{sec:comp}, we explain how the lists were computed
and present the theoretical background which assures that the computations were exhaustive. In Section~\ref{sec:doc}, information about the
format of the files available on \cite{online} is given.

We now proceed with a few comments on the census of $4$-HATs.
By a {\em vertex-stabiliser} of a vertex-transitive graph or digraph $\Gamma$, we mean the stabiliser of a vertex in $\Aut(\Gamma)$.
Even though it is known that a vertex-stabiliser of a $4$-HAT can be arbitrarily large (see \cite{DM05}),
not many examples of $4$-HATs with  vertex-stabilisers of order larger than $2$ were known, and all known examples had a very large number of vertices.
Recently, Conder and \v{S}parl (see also \cite{ConPotSpa}) discovered a $4$-HAT on $256$ vertices with  vertex-stabiliser of order $4$ and proved
that this is the smallest such example. This fact is confirmed by our census; in fact, the following theorem can be deduced from the census.

\begin{theorem}
\label{the:largeGv}
Amongst the 3246 $4$-HATs on at most $1000$ vertices, there are seventeen with vertex-stabiliser of order $4$, three with vertex-stabiliser of order $8$, and none with 
larger vertex-stabilisers. The smallest $4$-HAT with vertex-stabiliser of order $4$ has order $256$ and the smallest two with vertex-stabilisers of order $8$
have $768$ vertices; the third $4$-HAT with vertex-stabiliser of order $8$ has $896$ vertices.
\end{theorem}

Another curiosity about $4$-HATs is that those with a non-abelian vertex-stabiliser tend to be very rare (at least amongst the ``small'' graphs).
The first known $4$-HAT with a non-abelian vertex-stabiliser was discovered by Conder and Maru\v{s}i\v{c} (see~\cite{ConderMarusic}) and has $10752$ vertices. Further examples of $4$-HATs with non-abelian vertex-stabilisers were discovered recently (see \cite{ConPotSpa}),
including one with a vertex-stabiliser of order $16$. However, the one on $10752$ vertices remains the smallest known example. Using our list, the following fact is easily checked.

\begin{theorem}
\label{the:HATnab}
Every $4$-HAT with a non-abelian vertex-stabiliser has more than $1000$ vertices.
\end{theorem}

In fact, there are strong indications that the graph on $10752$ vertices discovered by Conder and Maru\v{s}i\v{c} is the smallest $4$-HAT with  a non-abelian vertex-stabiliser.

We will call a $4$-HAT with a non-solvable automorphism group a {\em non-solvable $4$-HAT}.
The first known non-solvable $4$-HAT was constructed by Maru\v{s}i\v{c} and Xu \cite{MarXu}; and its order is
 $7!/2$. An infinite family of non-solvable $4$-HATs was constructed later by Malni\v{c} and Maru\v{s}i\v{c} \cite{MalMar99}. The smallest member of this family has an even larger order, namely $11!/2$. To the best of our knowledge, no smaller non-solvable
$4$-HATs was known prior to the construction of our census. Perhaps surprisingly, small examples of non-solvable $4$-HATs seem not to be too rare,
as can be checked from our census. (The terms {\em radius},  {\em attachment number}, {\em alter-exponent}, and {\em alter-perimeter}
are defined in Sections \ref{sec:alt} and \ref{sec:rad}.)

\begin{theorem}
There are thirty-two non-solvable $4$-HATs with at most $1000$ vertices. The smallest one, named HAT[480,44], has order $480$, girth $5$, radius $5$,
attachment number $2$, alter-exponent $2$, and alter-perimeter $1$. It is non-Cayley and non-bipartite. 
\end{theorem}

Let us now continue with a few comments on the census of $2$-ATDs. All the undefined notions mentioned in the theorems below 
are explained in Sections~\ref{sec:not}, \ref{sec:alt} and \ref{sec:rad}. It is not surprising that, apart from the generalised wreath digraphs (see Section~\ref{sec:wreath} for the definition),
very few of the $2$-ATDs on at most $1000$ vertices are $2$-arc-transitive. In fact, the following can be deduced from the census.

\begin{theorem}
Out of the 26457 $2$-ATDs on at most $1000$ vertices, 961 are generalised wreath digraphs. Of the remaining 25496, only 1199 are $2$-arc-transitive (the smallest having order $18$), 
only 255 are $3$-arc-transitive (the smallest having order $42$), only 61 are $4$-arc-transitive (the smallest having order $90$), and only 5 are $5$-arc-transitive
(the smallest two having order $640$); none of them is $6$-arc-transitive.
\end{theorem}

Note that the non-existence of a $6$-arc-transitive non-generalised-wreath $2$-ATD on at most $1000$ vertices follows from a more general result (see
Corollary~\ref{cor:genlost}).

Recall that there is no $4$-HAT on at most $1000$ vertices with a non-abelian vertex-stabiliser   (Theorem~\ref{the:HATnab}). Consequently
(see Section~\ref{sec:ATvsGHAT}), every $2$-ATD on at most $1000$ vertices with a non-abelian vertex-stabiliser 
has an arc-transitive underlying graph; and there are indeed such examples.
In fact, the following holds (see Section~\ref{sec:def} for the definition of {\em self-opposite}).

\begin{theorem}
There are precisely forty-five $2$-ATDs on at most $1000$ vertices with a non-abelian vertex-stabiliser. 
They are all self-opposite, at least $3$-arc-transitive, have non-solvable automorphism groups,
and radius $3$. The smallest of these digraphs has order $42$, and the smallest that is $4$-arc-transitive has order $90$. 
There are no $5$-arc-transitive $2$-ATDs with a non-abelian vertex-stabiliser and order at most $1000$.
\end{theorem}

If a $2$-ATD is self-opposite, then the isomorphism between the digraph and its opposite digraph is an automorphism of the underlying graph,
making the underlying graph arc-transitive. Hence, self-opposite $2$-ATDs always yield arc-transitive $4$-GHATs. However, 
the converse is not always true: there are $2$-ATDs that are not self-opposite, but have an arc-transitive underlying graph.
In this case, the index of the automorphism group of the $2$-ATD in the automorphism group of its underlying graph must be larger than $2$
(for otherwise the former would be normal in the latter and thus any automorphism of the underlying graph would either preserve
the arc-set of the digraph, or map it to the arc-set of the opposite digraph). It is perhaps surprising that there are not many small examples of such behaviour.

\begin{theorem}
There are precisely fifty-two $2$-ATDs on at most $1000$ vertices that are not self-opposite but have an arc-transitive underlying graph.
The smallest two have order $21$. None of these digraphs is $2$-arc-transitive. The index of the automorphism group of these digraphs
in the automorphism group of the underlying graphs is always $8$.
\end{theorem}

\section{Notation and definitions}
\label{sec:not}

\subsection{Digraphs and graphs}
\label{sec:def}

A \emph{digraph} is an ordered pair $(V,A)$ where $V$ is a finite non-empty set and $A\subseteq V \times V$ is a binary relation on $V$. We say that $(V,A)$ is \emph{asymmetric} if $A$ is asymmetric, and we say that $(V,A)$ is a \emph{graph} if $A$ is irreflexive and symmetric. If $\Gamma=(V,A)$ is a digraph, then
we shall refer to the set $V$ and the relation $A$ as the {\em vertex-set} and the {\em arc-set} of $\Gamma$, and denote them by
$\V(\Gamma)$ and $\A(\Gamma)$, respectively. Members of $V$ and $A$ are called {\em vertices} and {\em arcs}, respectively. If $(u,v)$ is an arc of a digraph $\Gamma$, then $u$ is called the {\em tail}, and $v$ the {\em head} of $(u,v)$. If $\Gamma$ is a graph, then the unordered pair $\{u,v\}$ is called an
{\em edge} of $\Gamma$ and the set of all edges of $\Gamma$ is denoted $\E(\Gamma)$.

If $\Gamma$ is a digraph, then the {\em opposite digraph} $\oppG$ has  vertex-set $\V(\Gamma)$ and arc-set $\{(v,u) :  (u,v) \in \A(\Gamma)\}$. The {\em underlying graph} of $\Gamma$ is the graph with vertex-set $\V(\Gamma)$ and with arc-set $\A(\Gamma) \cup \A(\oppG)$.
A digraph is called {\em connected} provided that its underlying graph is connected.

Let $v$ be a vertex of a digraph $\Gamma$. Then the {\em out-neighbourhood} of $v$ in $\Gamma$, denoted by $\Gamma^+(v)$,
is the set of all vertices $u$ of $\Gamma$ such that $(v,u) \in \A(\Gamma)$, and similarly, the {\em in-neighbourhood} $\Gamma^-(v)$ is defined 
as the set of all vertices $u$ of $\Gamma$ such that $(u,v) \in \A(\Gamma)$.
Further, we let $\val^+(v) = |\Gamma^+(v)|$ and $\val^-(v) = |\Gamma^-(v)|$ be the {\em out-valence} and {\em in-valence} of $\Gamma$, respectively.
If there exists an integer $r$ such that $\val^+(v) = \val^-(v) = r$ for every $v\in \V(\Gamma)$, then we say that $\Gamma$ is {\em regular} of {\em valence} $r$, or
simply that $\Gamma$ is an {\em $r$-valent} digraph.

An $s$-arc of a digraph $\Gamma$ is an $(s+1)$-tuple $(v_0,v_1, \ldots, v_s)$ of vertices of $\Gamma$, such that $(v_{i-1},v_i)$
is an arc of $\Gamma$ for every $i\in \{1,\ldots,s\}$ and
$v_{i-1}\not=v_{i+1}$ for every $i\in \{1,\ldots,s-1\}$.
If $x=(v_0,v_1, \ldots, v_s)$ is an $s$-arc of  $\Gamma$, then every $s$-arc of the form $(v_1, v_2, \ldots, v_s,w)$ is called a {\em successor} of $x$.

An \emph{automorphism} of a digraph $\Gamma$ is a permutation of $\V(\Gamma)$ which preserves
the arc-set $\A(\Gamma)$. 
Let $G$ be a subgroup of the full automorphism group $\Aut(\Gamma)$ of $\Gamma$.
We say that $\Gamma$ is \emph{$G$-vertex-transitive} or \emph{$G$-arc-transitive} provided that
$G$ acts transitively on $\V(\Gamma)$ or $\A(\Gamma)$, respectively. Similarly, we say
that $\Gamma$ is \emph{$(G,s)$-arc-transitive} if $G$ acts
transitively on the set of $s$-arcs of $\Gamma$. If $\Gamma$ is a graph, we say that it is \emph{$G$-edge-transitive} provided that $G$ acts transitively on $\E(\Gamma)$.
When $G=\Aut(\Gamma)$, the prefix $G$ in the above notations is usually omitted.

If $\Gamma$ is a digraph and $v\in \V(\Gamma)$, then a {\em $v$-shunt} is an automorphism of $\Gamma$ which maps $v$ to an out-neighbour of $v$.

\subsection{From $4$-GHATs to $2$-ATDs and back}
\label{sec:ATvsGHAT}

If $\Gamma$ is a connected $4$-valent  $(G,\frac{1}{2})$-arc-transitive graph, then $G$ has two orbits on the arc-set of $\Gamma$, opposite to each other,
each orbit having the property that each vertex of $\Gamma$ is the head of precisely two arcs, and also the tail of
precisely two arcs of the orbit. By taking any of these two orbits as an arc-set of a digraph on the same vertex-set, one thus obtains a 2-ATD whose underlying graph is $\Gamma$, and admitting $G$ as an arc-transitive group of automorphisms.

Conversely, the underlying graph of a $G$-arc-transitive $2$-ATD is a $(G,\frac{1}{2})$-arc-transitive $4$-GHAT. In this sense the study of $4$-GHATs is  equivalent to the study of $2$-ATDs. 

In Section~\ref{sec:comp}, we explain how a complete list of all $2$-ATDs on at most $1000$ vertices
was obtained. The above discussion shows how this yields a complete list of all $4$-GHATs on at most $1000$ vertices.

\subsection{Generalised wreath digraphs}
\label{sec:wreath}

Let $n$ be an integer with $n\ge 3$, let $V=\ZZ_n\times \ZZ_2$, and let $A=\{ ((i,a),(i+1,b)) : i \in \ZZ_n, a,b\in \ZZ_2\}$.
The asymmetric digraph $(V,A)$ is called a {\em wreath digraph}
and denoted by $\vW_n$. 

If $\Gamma$ is a digraph and $r$ is a positive integer, then the {\em $r$-th partial line digraph} of $\Gamma$, denoted $\Pl^r(\Gamma)$, 
is the digraph with vertex-set equal to the set of $r$-arcs of $\Gamma$ and with
$(x,y)$ being an arc of $\Pl^r(\Gamma)$ whenever $y$ is a successor of $x$. If $r=0$, then we let $\Pl^r(\Gamma)=\Gamma$.

Let $r$ be a positive integer. 
The $(r-1)$-th partial line digraph $\Pl^{r-1}(\vW_n)$ of the wreath digraph $\vW_n$ is denoted by $\vW(n,r)$ and called a {\em generalised wreath digraph}.
Generalised wreath digraphs were first introduced in \cite{PraHATD}, where $\vW(n,r)$ was denoted $C_n(2,r)$.
It was proved there that $\Aut(\vW(n,r)) \cong C_2\wr C_n$  and that $\Aut(\vW(n,r))$  acts transitively on the $(n-r)$-arcs but not on the
$(n-r+1)$-arcs of $\vW(n,r)$~\cite[Theorem 2.8]{PraHATD}. In particular, $\vW(n,r)$ is arc-transitive if and only if $n\ge r+1$.
Note that $|\V(\vW(n,r))| = n2^{r}$, and thus $|\Aut(\vW(n,r))_v| = n2^n/n2^{r} = 2^{n-r}$.

The underlying graph of a generalised wreath digraph will be called a {\em generalised wreath graph}.

\subsection{Coset digraphs}

Let $G$ be a group generated by a core-free subgroup $H$ and an element $g$ with $g^{-1} \not\in HgH$.
One can construct the {\em coset digraph}, denoted $\Cos(G,H,g)$,
whose vertex-set is the set $G/H$ of right cosets of $H$ in $G$, and where $(Hx,Hy)$ is an arc if and only if $yx^{-1} \in HgH$.
Note that the condition $g^{-1} \not\in HgH$ guarantees that the arc-set is an asymmetric relation. Moreover, since $G=\langle H,g\rangle$,
the digraph $\Cos(G,H,g)$ is connected.

The digraph $\Cos(G,H,g)$ is $G$-arc-transitive (with $G$ acting upon $G/H$ by right multiplication),
and hence $\Cos(G,H,g)$ is a $G$-arc-transitive and $G$-vertex-transitive digraph with $g$ being a $v$-shunt.
On the other hand, it is folklore that every such graph arises as a coset digraph.

\begin{lemma}
\label{lem:coset}
If $\Gamma$ is a connected $G$-arc-transitive and $G$-vertex-transitive
digraph, $v$ is a vertex of $\Gamma$, and $g$ is a $v$-shunt contained in $G$, then $\Gamma \cong \Cos(G,G_v,g)$. 
\end{lemma}

\section{Constructing the census}
\label{sec:comp}

If $\Gamma$ is a $G$-vertex-transitive digraph with $n$ vertices, then $|G| = n|G_v|$.
If one wants to use the coset digraph construction to obtain all $2$-ATDs on $n$ vertices,
one thus needs to consider all groups $G$ of order $n|G_v|$ that can act as arc-transitive groups of
$2$-ATDs.
In order for this approach to be practical, two issues must be resolved:

First, one must get some  control over $|G|$ and thus over $|G_v|$. (Recall that in $\vW(n,r)$, $|G_v|$ can grow exponentially with $|\V(\vW(n,r))|$, as $n\to \infty$ and $r$ is fixed). Second, one must obtain enough structural information about $G$ to be able to construct all possibilities.

Fortunately, both of these issues were resolved successfully. The problem of bounding $|G_v|$ was resolved in a recent paper \cite{genlost} and details can be found in Section~\ref{sec:bound}. The second problem was dealt with in \cite{MarNed3}, and later, in greater generality in \cite{PotVer} (both of these papers rely heavily on a group-theoretical result of Glauberman \cite{Glaub2}); the summary of relevant results is given in Section~\ref{sec:type}.

\subsection{Bounding the order of the vertex-stabiliser}
\label{sec:bound}

The crucial result that made our compilation of a complete census of all small $2$-ATDs possible is Theorem~\ref{thm:genlost}, stated below,
which shows that the generalised wreath digraphs (defined in Section~\ref{sec:wreath})
are very special in the sense of having large vertex-stabilisers. In fact, together with the correspondence described in Section~\ref{sec:ATvsGHAT}, \cite[Theorem~9.1]{genlost} has the following corollary:

\begin{theorem}
\label{thm:genlost}
Let $\Gamma$ be a $G$-arc-transitive $2$-ATD on at most $m$ vertices and let $t$ be the largest integer such that $m> t 2^{t+2}$.
Then  one of the following occurs:
\begin{enumerate}
\item $\Gamma\cong \vW(n,r)$ for some $n\ge 3$ and $1\le r \le n-1$,
\item $|G_v| \le \max\{16,2^t\}$,
\item $(\Gamma,G)$ appears in the last line of \cite[Table~5]{genlost}. In particular, $|\V\Gamma|=8100$.
\end{enumerate}
\end{theorem}

The following is an easy corollary:

\begin{corollary}
\label{cor:genlost}
Let $\Gamma$ be a $G$-arc-transitive $2$-ATD on at most $1000$ vertices. Then either $|G_v| \le 32$ or $\Gamma\cong \vW(n,r)$ for some $n\ge 3$ and $1\le r \le n-1$.
\end{corollary}

\subsection{Structure of the vertex-stabiliser}
\label{sec:type}

\begin{definition}\label{defdef}
Let $s$ and $\alpha$ be positive integers satisfying $\frac{2}{3}s\le \alpha \le s$,
and let $c$ be a function assigning a value $c_{i,j}\in \{0,1\}$ to each pair of integers $i,j$ with $\alpha \le j \le s-1$ and $1\le i \le 2\alpha-2s+j+1$. Let $A_{s,\alpha}^c$ be the group generated by  $\{x_0, x_1, \ldots, x_{s-1}, g\}$ and subject to the defining relations:
\begin{itemize}
\item $x_0^2 = x_1^2 = \cdots = x_{s-1}^2 = 1$;
\item $x_i^g=x_{i+1}$ for $i\in\{0,1,\ldots,s-2\}$;
\item if $j < \alpha$, then $[x_0,x_j] = 1$;
\item if $j\ge \alpha$, then $[x_0,x_j] =  x_{s-\alpha}^{c_{1,j}}\, x_{s-\alpha+1}^{c_{2,j}}\,\cdots\, x_{j-s+\alpha}^{c_{2\alpha-2s+j+1,j}}$.
\end{itemize}
\end{definition}
Furthermore, let $\cA_{s,\alpha}$ be the family of all groups $A_{s,\alpha}^c$ for some $c$. It was proved in \cite{MarNed3} (see also \cite{PotVer}) that every group $G$ acting arc-transitively 
on a $2$-ATD is isomorphic to a quotient of some $A_{s,\alpha}^c$. More precisely, the following can be deduced from \cite{MarNed3} or \cite{PotVer}.

\begin{theorem}
\label{thm:structure}
Let $\Gamma$ be a $G$-arc-transitive $2$-ATD, let $v\in \V(\Gamma)$ and let $s$ be the largest integer such that $G$ acts transitively
on the set of $s$-arcs of $\Gamma$. Then there exists an integer $\alpha$ satisfying $\frac{2}{3}s\le \alpha \le s$, a function $c$ as in Definition~\ref{defdef},
and an epimorphism $\wp\colon A_{s,\alpha}^c \to G$, which maps the group $\langle x_0, \ldots, x_{s-1}\rangle$ isomorphically onto $G_v$ 
and the generator $g$ to some $v$-shunt in $G$. In particular, $|G_v|=2^s$.
\end{theorem}

In this case, we will say that $(\Gamma,G)$ is of {\em type} $A_{s,\alpha}^c$, and call the group $A_{s,\alpha}^c$ the {\em universal group} of
the pair $(\Gamma,G)$.

For $s$, $\alpha$, and a function $c$ satisfying the conditions of Definition~\ref{defdef}, let
$c'$ be the function defined by $c'_{i,j} = c_{2\alpha-2s+j+2\,-\,i,j}$. The relationship between $c$ and $c'$ can
be visualised as follows: if one fixes the index $j$ and views the function $i\mapsto c_{i,j}$ as the sequence
$[c_{1,j}, c_{2,j}, \ldots, c_{2\alpha-2s+j+1,j}]$, then the sequence for $c'$ is obtained by reversing the one for $c$.
If $\tG = A_{s,\alpha}^c$ then we denote the {\em reverse type} $A_{s,\alpha}^{c'}$ by $\tG^{{\rm opp}}$.
 
Observe that if $(\Gamma,G)$ is of type $\tG$, then $(\oppG,G)$ is of type $\tG^{{\rm opp}}$.
A class of groups, obtained from $\cA_{s,\alpha}$ by taking only one group in each pair $\{\tG,\tG^{{\rm opp}}\}$, $\tG\in \cA_{s,\alpha}$, will
be denoted $\cA_{s,\alpha}^{{\rm red}}$. (Note that some groups $\tG$ might have the property that $\tG = \tG^{{\rm opp}}$.)

In view of Corollary~\ref{cor:genlost}, we shall be mainly interested in the universal groups $A_{s,\alpha}^c$ with $s\le 5$ (as, excluding generalised wreath digraphs, these are the only types of $2$-ATDs of order at most $1000$).
We list the relevant classes $\cA_{s,\alpha}^{\rm{red}}$ for $s\le 5$ explicitly in
Table~\ref{tab:1}. Groups in $\cA_{s,\alpha}^{\rm{red}}$, for a fixed $s$
 will be named by $\am_s^i$, where $i$ will be a positive integer, where groups with larger $\alpha$ will be indexed with lower $i$.
 Also, the generators $x_0, x_1, x_2, x_3$, and $x_4$ will be denoted $a$, $b$, $c$, $d$, and $e$, respectively.

\begin{table}[hhh]
\begin{center}
\begin{small}
\begin{tabular}{|c|c|}
\hline
\phantom{$\overline{\overline{G_j^G}}$}
name & $\tilde{G}$ \\
\hline\hline
$\am_1^1$ &
\begin{tabular}{l}
\phantom{$\overline{\overline{G_j^G}}$}
$\la a,g \mid a^2 \ra$ \\
\end{tabular} 
\\
\hline
$\am_2^1$ &
\begin{tabular}{l}
\phantom{$\overline{\overline{G_j^G}}$}
$\la a,b,g \mid a^2,b^2,a^gb, [a,b] \ra$ \\
\end{tabular} 
\\
\hline
$\am_3^1$ &
\begin{tabular}{l}
\phantom{$\overline{\overline{G_j^G}}$}
$\la a,b,c,g \mid a^2,b^2,c^2,a^gb,b^gc,[a,b],[a,c] \ra$ \\
\end{tabular} 
\\
\hline
$\am_3^2$ &
\begin{tabular}{l}
\phantom{$\overline{\overline{G_j^G}}$}
$\la a,b,c,g \mid a^2,b^2,c^2,a^gb,b^gc,[a,b],[a,c]b \ra$ \\
\end{tabular} 
\\
\hline
$\am_4^1$ &
\begin{tabular}{l}
\phantom{$\overline{\overline{G_j^G}}$}
$\la a,b,c,d,g \mid a^2,b^2,c^2,d^2,a^gb, b^gc, c^gd, [a,b],[a,c],[a,d] \ra$ \\
\end{tabular} 
\\
\hline
$\am_4^2$ &
\begin{tabular}{l}
\phantom{$\overline{\overline{G_j^G}}$}
$\la a,b,c,d,g \mid a^2,b^2,c^2,d^2,a^gb, b^gc, c^gd, [a,b],[a,c],[a,d]b \ra$ \\
\end{tabular} 
\\
\hline
$\am_4^3$ &
\begin{tabular}{l}
\phantom{$\overline{\overline{G_j^G}}$}
$\la a,b,c,d,g \mid a^2,b^2,c^2,d^2,a^gb, b^gc, c^gd, [a,b],[a,c],[a,d]bc \ra$ \\
\end{tabular} 
\\
\hline
$\am_5^1$ &
\begin{tabular}{l}
\phantom{$\overline{\overline{G_j^G}}$}
$\la a,b,c,d,e,g \mid a^2,b^2,c^2,d^2,e^2,d^2,a^gb,b^gc,c^gd,d^ge,[a,b],[a,c],[a,d],[a,e] \ra$ \\
\end{tabular} 
\\
\hline
$\am_5^2$ &
\begin{tabular}{l}
\phantom{$\overline{\overline{G_j^G}}$}
$\la a,b,c,d,e,g \mid a^2,b^2,c^2,d^2,e^2,d^2,a^gb,b^gc,c^gd,d^ge,[a,b],[a,c],[a,d],[a,e]b \ra$ \\
\end{tabular} 
\\
\hline
$\am_5^3$ &
\begin{tabular}{l}
\phantom{$\overline{\overline{G_j^G}}$}
$\la a,b,c,d,e,g \mid a^2,b^2,c^2,d^2,e^2,d^2,a^gb,b^gc,c^gd,d^ge,[a,b],[a,c],[a,d],[a,e]c \ra$ \\
\end{tabular} 
\\
\hline
$\am_5^4$ &
\begin{tabular}{l}
\phantom{$\overline{\overline{G_j^G}}$}
$\la a,b,c,d,e,g \mid a^2,b^2,c^2,d^2,e^2,d^2,a^gb,b^gc,c^gd,d^ge,[a,b],[a,c],[a,d],[a,e]bc \ra$ \\
\end{tabular} 
\\
\hline
$\am_5^5$ &
\begin{tabular}{l}
\phantom{$\overline{\overline{G_j^G}}$}
$\la a,b,c,d,e,g \mid a^2,b^2,c^2,d^2,e^2,d^2,a^gb,b^gc,c^gd,d^ge,[a,b],[a,c],[a,d],[a,e]bd \ra$ \\
\end{tabular} 
\\
\hline
$\am_5^6$ &
\begin{tabular}{l}
\phantom{$\overline{\overline{G_j^G}}$}
$\la a,b,c,d,e,g \mid a^2,b^2,c^2,d^2,e^2,d^2,a^gb,b^gc,c^gd,d^ge,[a,b],[a,c],[a,d],[a,e]bcd \ra$ \\
\end{tabular} 
\\
\hline                  
                
\end{tabular}
\end{small}
\caption{Universal groups of $2$-ATDs with $|\tilde{G}_v| \le 32$}     
\label{tab:1}
\end{center}
\end{table}

\subsection{The algorithm and its implementation}

We now have all the tools required to present a practical algorithm  that takes an integer $m$ as input and returns 
a complete list of all $2$-ATDs on at most $m$ vertices (see Algorithm~1).
It is based on the fact that every such digraph can be obtained as
a coset digraph of some group $G$ (see Lemma~\ref{lem:coset}), and that $G$ is in fact an epimorphic image of some
group $A_{s,\alpha}^c$ (see Theorem~\ref{thm:structure}) with $G_v$ and the shunt being the corresponding 
images of $\langle x_0, \ldots, x_{s-1}\rangle$ and $g$ in $A_{s,\alpha}^c$. 

Moreover, if  $\Gamma$ is not a generalised wreath digraph
or the exceptional digraph on $8100$ vertices mentioned in part 3 of Theorem~\ref{thm:genlost}, then the parameter $s$ satisfies
$s2^{s+2}<m$, and the order of the epimorphic image $G$ is bounded by $2^s m$
(see Theorem~\ref{thm:genlost}). The algorithm thus basically boils down to the task of finding normal subgroups of bounded index in
the finitely presented groups $A_{s,\alpha}^c$.

\begin{algorithm}
\begin{algorithmic}
\label{alg:main}
\REQUIRE positive integer $m$
\ENSURE $\D = \{ \Gamma : \Gamma$ is 2-ATD, $|V(\Gamma)| \le m\}$
\STATE Let $t$ be the largest integer such that $m> t 2^{t+2}$;
\STATE Let $\D$ be the list of all arc-transitive generalised wreath digraphs on at most $m$ vertices;
\STATE If $m\ge 8100$, add to $\D$ the exceptional digraph $\Gamma$ on $8100$ vertices, mentioned in part 3 of Theorem~\ref{thm:genlost};
\FOR{$s\in\{1,\ldots,\max\{4,t\}\}$}
\FOR{$\alpha \in \{ \lceil\frac{2}{3}s\rceil, \lceil\frac{2}{3}s\rceil+1, \ldots, s \}$}
\FOR{$\tG \in \cA_{s,\alpha}^{{\rm red}}$}
\STATE Let $\N$ be the set of all normal subgroups of $\tG$ of index at most $2^sm$;
\FOR{$N \in \N$}
\STATE Let  $G:=\tG/N$ and let $\wp \colon \tG \to G$ be the quotient projection;
\STATE Let $H:=\wp(\langle x_0, \ldots, x_{s-1}\rangle)$;
\IF{$H$ is core-free in $G$ \AND $|H| = 2^s$ \AND $\wp(g)^{-1} \not \in H\wp(g)H$}
  \STATE Let $C :=\cos(G,H,\wp)$;
  \FOR{$\Gamma\in \{C,C^{\rm{opp}}\}$}
  \IF{$\Gamma$ is not isomorphism to any of the digraphs in $\D$}
     \STATE add $\Gamma$  to the list $\D$;
  \ENDIF   
  \ENDFOR
\ENDIF
\ENDFOR
\ENDFOR
\ENDFOR
\ENDFOR
\end{algorithmic}\caption{~$2$-ATDs on at most $m$ vertices.}
\end{algorithm}

Practical implementations of this algorithm have several limitations. 
First, the best known algorithm for finding normal subgroups of low index in
a finitely presented group is an algorithm due to  Firth and  Holt~\cite{Firth}. The only publicly available implementation is the {\tt LowIndexNormalSubgroups} routine in {\sc Magma} \cite{Magma} and the most recent version allows one to compute only the normal subgroups of index at most $5\cdot 10^5$; hence only automorphisms groups of order $5\cdot 10^5$ can possibly be
obtained in this way.

More importantly, even when only normal subgroups of relatively small index need to be computed, some finitely presented groups are computationally difficult.
For example, finding all normal subgroups of index at most $2048$ of
the group $A_1^1\cong C_2*C_\infty$ seems to represent a considerable challenge for the {\tt LowIndexNormalSubgroups} routine in {\sc Magma}.
In order to overcome this problem, we have used a recently computed catalogue of all $(2,*)$-groups of order at most $6000$ \cite{2star}, 
where by a {\em $(2,*)$-group} we mean any group generated by an involution $x$ and one other element $g$. Since $A_1^1$ is a $(2,*)$-group and every non-cyclic quotient of  a $(2,*)$-group is also a $(2,*)$-group, this catalogue can be used to obtain all the quotients of $A_1^1$ of order up to $6000$.
Consequently, all $2$-ATDs admitting an arc-regular group of automorphisms of order at most $3000$ can be obtained.
Similarly, since $A_2^1$ is also a $(2,*)$-group, we can use this catalogue to obtain all the $2$-ATDs of order at most $1500$
admitting an arc-transitive group $G$ with $|G_v|=4$.

Like $A_1^1$ and $A_2^1$, the groups with $\langle x_0, \ldots, x_{s-1}\rangle$ abelian (namely those with $\alpha=s$ and $c_{i,j}=0$ for all $i,j$)  
are also computationally  very difficult.
One can make the task easier by dividing it into cases, where the order of $g$ is fixed in each case. Since $g$ represents a shunt, it can be
proved that its order cannot exceed the order of the digraph (see, for example, \cite[Lemma 13]{cubiccensus}). Cases can then be run in parallel on a multi-core algorithm.

\section{The census and accompanying data}
\label{sec:doc}

Using Algorithm~1, we found that there are exactly  26457  $2$-ATDs of order up to $1000$. Following the recipe explained in Section~\ref{sec:ATvsGHAT},
  we have also computed all the $4$-GHATs, which we split in two lists: $4$-HATs and arc-transitive $4$-GHATs.

The data about these graphs, together with {\sc Magma} code
that generates them, is available on-line at \cite{online}. The package contains ten files. The file
``Census-ATD-1k-README.txt'' is a text file containing information similar to the information in this section. The remaining nine files
come in groups of three, one group for each of the three lists ($2$-ATDs, arc-transitive $4$-GHATs, $4$-HATs). In each groups,
there is a $*$.mgm file, a $*$.txt file and a $*$.csv file. 

The $*$.mgm file contains {\sc Magma} code that generates
the corresponding digraphs. After loading the file in {\sc Magma}, a double sequence is generated (named either ATD, GHAT, or HAT, depending on the file). The length of each double sequence is $1000$ and
the $n$-th component of the sequence is the sequence of all the corresponding digraphs of order $n$, with the exception of the
generalised wreath digraphs.
Thus, ATD[32,2] will return the second of the four non-generalised-wreath $2$-ATDs on $32$ vertices 
(the ordering of the digraphs in the sequence ATD[32] is arbitrary).
In order to include the generalised wreath digraphs into the corresponding sequence, one can call the procedure {\tt AddGWD($\sim$ATD,GWD)}
in the case of the $2$-ATDs, or {\tt AddGWG($\sim$GHAT,GWG)} in the case of the $4$-GHATs (note that a generalised wreath graph is never \half-arc-transitive).

The $*$.txt file contains the list of neighbours of each digraph. This file is needed when the $*$.mgm file is loaded into {\sc Magma}, but, being an ASCII file,
it can be used also by other computer systems to reconstruct the digraphs. For the details of the format, see the ``README'' file.

Finally, the $*$.csv file is a ``comma separated values'' file representing a spreadsheet containing some precomputed 
graph invariants. We shall first introduce some of these invariants and then
discuss each $*$.csv separately.

\subsection{Walks and cycles}

Let $\Gamma$ be a digraph. 
 A {\em walk} of length $n$ in $\Gamma$ is an $(n+1)$-tuple
$(v_0,v_1,\ldots,v_n)$ of vertices of $\Gamma$ such that, for any $i\in\{1,\ldots n\}$, either
$(v_{i-1},v_i)$ or $(v_i,v_{i-1})$ is an arc of $\Gamma$. The walk is {\em closed} if $v_0=v_n$ and {\em simple}
if the vertices $v_i$ are pairwise distinct (with the possible exception of the first and the last vertex when the walk is closed).

A closed simple walk in $\Gamma$ is called a {\em cyclet}.
The {\em inverse} of a cyclet $(v_0, \ldots, v_{n-1},v_0)$ is the cyclet $(v_0, v_{n-1}, \ldots, v_1,v_0)$, and 
a cyclet $(v_0, \ldots, v_{n-1},v_0)$ is said to be
a {\em shift} of a cyclet $(u_0,\ldots,u_{n-1},u_0)$ provided that there exists $k\in \ZZ_n$ such that $u_i = v_{i+k}$ for all $i\in \ZZ_n$.
Two cyclets $W$ and $U$ are said to be {\em congruent} provided that $W$ is a shift of either $U$ or the inverse of $U$.
The relations of ``being a shift of'' and ``being congruent to'' are clearly equivalence relations, and their equivalence classes are called {\em oriented cycles}
and {\em cycles}, respectively. With a slight abuse of terminology, we shall sometimes identify a (oriented) cycle with any of its representatives.

\subsection{Alter-equivalence,  alter-exponent, alter-perimeter, and alter-sequence}
\label{sec:alt}
Let $\Gamma$ be an asymmetric digraph. The {\em signature} of a walk $W=(v_0,v_1,\ldots,v_n)$ is
an $n$-tuple $(\epsilon_1, \epsilon_2, \ldots, \epsilon_n)$, where $\epsilon_i=1$ if $(v_{i-1},v_i)$ is an arc of
$\Gamma$, and $\epsilon=-1$ otherwise. The signature of a walk $W$ will be denoted by $\sigma(W)$.
The sum of all the integers in $\sigma(W)$ is called the {\em sum} of the walk $W$ and denoted by $s(W)$; similarly,
the $k^{th}$ {\em partial sum $s_k(W)$}  is the sum of the initial walk
$(v_0,v_1,\ldots,v_k)$ of length $k$. By convention, we let $s_0(W)=0$. 

The {\em tolerance} of a walk $W$ of length $n$, denoted $T(W)$, is the set $\{s_k(W) : k\in \{0,1,\ldots, n\}\}$.
Observe that the tolerance of a walk is always an interval of integers containing $0$.
Let $t$ be a positive integer or $\infty$. We say that two vertices $u$ and $v$ of $\Gamma$ are {\em alter-equivalent with
tolerance $t$} if there is a walk from $u$ to $v$
with sum $0$ and tolerance contained in $[0,t]$; we shall then write $u \alt_t v$. 
The equivalence class of $\alt_t$ containing a vertex $v$ will be denoted by $\alt_t(v)$.

Since we assume that $\Gamma$ is a finite digraph, there exists
an integer $e\geq 0$ such that $\alt_e = \alt_{e+1}$ (and then $\alt_e = \alt_\infty$).
The smallest such integer is called the {\em alter-exponent} of 
$\Gamma$ and denoted by $\exp(\Gamma)$.

The number of equivalence classes of $\alt_{\infty}$ is called the {\em alter-perimeter} of $\Gamma$.
The name originates from the fact that the quotient digraph of $\Gamma$ with respect to $\alt_\infty$ is
either a directed cycle or the complete graph $K_2$ or the graph $K_1$ with one vertex.

If $e$ is the alter-exponent of a (vertex-transitive) digraph $\Gamma$, then the finite sequence $[|\alt_1(v)|, |\alt_2(v)|, \ldots, |\alt_e(v)|]$ is called
the {\em alter-sequence} of $\Gamma$.

Several interesting properties of the alter-exponent can be proved (see \cite{bridge} for example).
For example, if $\Gamma$ is connected and $G$-vertex-transitive, then 
$\exp(\Gamma)$ is the smallest positive integer $e$ such that
the setwise stabiliser $G_{\alt_e(v)}$ is normal in $G$.
The group $G_{\alt_e(v)}$ is the group generated by all
vertex-stabilisers in $G$ and $G/G_{\alt_e(v)}$ is a cyclic group.

All notions defined in this section for digraphs generalise to half-arc-transitive graphs, where instead of the graph
one of the two natural arc-transitive digraphs are considered. As was shown in \cite{bridge},  all the parameters defined here
remain the same if instead of a digraph, its opposite digraph is considered.

\subsection{Alternating cycles -- radius and attachment number}
\label{sec:rad}
A walk $W$ in an asymmetric digraph is called {\em alternating} if its tolerance is either $[0,1]$ or $[-1,0]$ (that is, if the signs in its signature alternate).
Similarly, a cycle is called {\em alternating} provided that any (and thus every) of its representatives is an alternating walk.

This notion was introduced in \cite{Mar98} and used to classify the so-called
{\em tightly attached} $4$-GHATs. The concept of alternating cycles was explored further in a number of papers on
$4$-HATs (see for example \cite{MarPra,Spa08}).

Let $\Gamma$ be a $2$-ATD, let $\C$ be the set of all alternating cycles of $\Gamma$, and let $G=\Aut(\Gamma)$.
The set $\C$ is clearly preserved by the action of $G$ upon the cycles of $\Gamma$. Moreover, since $\Gamma$ is arc-transitive,
$G$ acts transitively on $\C$. In particular, all the alternating cycles of $\Gamma$ are of equal length. Half of the length of an alternating cycle is
called the {\em radius} of $\Gamma$.

Since  $\Gamma$ is $2$-valent, every vertex of $\Gamma$ belongs to precisely two alternating cycles. It thus follows from vertex-transitivity of $\Gamma$
that any (unordered) pair of intersecting cycles can be mapped to any other such pair, implying that there exists a constant $a$ such that any two cycles
meet either in $0$ or in $a$ vertices. The parameter $a$ is then called the {\em attachment number} of $\Gamma$.
In general, the attachment number divides the length of the alternating cycle (twice the radius), and
there are digraphs where $a$ equals this length; they were classified in
\cite[Proposition 2.4]{Mar98}, where it was shown that their underlying graphs are always arc-transitive. 
A $2$-valent asymmetric digraph with attachment number $a$ is called {\em tightly attached} if $a$ equals the radius, is called {\em antipodally attached} if $a=2$,
and is called {\em loosely attached} if $a=1$. Note that tightly attached $2$-ATDs are precisely those with alter-exponent 1.

\subsection{Consistent cycles}

Let $\Gamma$ be a graph and let $G\le \Aut(\Gamma)$.
A (oriented) cycle $C$ in a graph $\Gamma$ is called $G$-consistent provided that there exists $g\in G$ that preserves $C$ and acts upon it as a $1$-step rotation.
A $G$-orbit of $G$-consistent oriented cycles is said to be {\em symmetric} if it contains the inverse of any (and thus each) of its members, and is {\em chiral} otherwise.

Consistent oriented cycles were first introduced by Conway in a public lecture \cite{con} (see also \cite{biggs,mik,overlap}). Conway's original result states that in an
arc-transitive graph of valence $d$, the automorphism group of the graph has exactly $d-1$ orbits on the set of oriented cycles. 
In particular, if $\Gamma$ is $4$-valent and $G$-arc-transitive, then there are precisely three $G$-orbits of $G$-consistent oriented cycles. Since chiral orbits of $G$-consistent
cycles  come in pairs of mutually inverse oriented cycles, this implies that there must be at least one symmetric orbit, while the other two are either both chiral or both symmetric.

Conway's result was generalised in \cite{ByK} to the case of \half-arc-transitive graphs by showing that if $\Gamma$ is a $4$-valent $(G,\frac{1}{2})$-arc-transitive graph,
then there are precisely four $G$-orbits of $G$-consistent oriented cycles, all of them chiral. These four orbits of oriented cycles thus constitute precisely two
$G$-orbits of $G$-consistent (non-oriented) cycles.

\subsection{Metacirculants}

A {\em metacirculant} is a graph whose automorphism group contains a vertex-transitive metacyclic group $G$, generated by $\rho$ and $\sigma$, such that
the cyclic group $\langle \rho \rangle$ is semiregular on the vertex-set of the graph, and is normal in $G$.
Metacirculants were first defined by Alspach and Parsons \cite{AlsPar}, and metacirculants admitting 
\half-arc-transitive groups of automorphisms were first investigated in \cite{sajna}. Recently, the interesting problem of classifying all $4$-HATs that are metacirculants
was considered in \cite{MarSpa,Spa09,Spa10}. Such $4$-HATs fall into four (not necessarily disjoint) classes (called Class I, Class II, Class III, and Class IV),
 depending on the structure of the quotient by the orbits of the semiregular element $\rho$. 
 For a precise definition of the {\em class} of a $4$-HAT metacirculant see, for example, \cite[Section 2]{Spa09}.
 Since a given $4$-HAT may admit several vertex-transitive metacyclic groups, a fixed graph can fall into several of these four classes.
 Several interesting facts about $4$-HAT metacirculants are known. For example, tightly attached $4$-HATs are precisely the $4$-HATs that are metacirculants of Class I.
 
\subsection{The data on $2$-ATDs}
\label{sec4.6}

The ``Census-ATD-1k-data.csv'' file concerns $2$-ATDs. Each line of the file represents one of the digraphs in the census, and
has 19 fields described below. Since this file is in ``csv'' format, every occurrence of a comma in a field is substitute with a semicolon. 

\begin{itemize}
\item {\tt Name}: the name of the digraph (for example, ATD[32,2]);
\item {\tt $|$V$|$}: the order of the digraph;
\item {\tt SelfOpp}:	contains ``yes" if the digraph is isomorphic to its opposite digraph and ``no" otherwise;
\item {\tt Opp}:		the name of the opposite digraph (the same as ``Name" if the digraph is self-opposite);
\item {\tt IsUndAT}:	``yes" if the underlying graph is arc-transitive, ``no'' otherwise;
\item {\tt UndGrph}:	the name of the underlying graph, as given in the files ``Census-HAT-1k-data.csv'' and ``Census-GHAT-1k-data.csv'' --
		if the underlying graph is generalized wreath, then this is indicated by, say, ``GWD(m,k)" where $m$ and $k$ are the defining parameters.
\item {\tt s}:	the largest integer $s$, such that the digraph is $s$-arc-transitive;
\item {\tt GvAb}:	``Ab'' if the vertex-stabiliser in the automorphism group of the digraph is abelian, otherwise ``n-Ab'';
\item {\tt $|$Tv:Gv$|$}:  the index of  the automorphism group $G$ of the digraph in the smallest arc-transitive group $T$ of the underlying graph that contains $G$ -- if there's no such group $T$, then $0$;
\item {\tt $|$Av:Gv$|$}: the index of the automorphism group of the digraph in the automorphism group of the underlying graph;
\item {\tt Solv}:		this field contains ``solve" if the automorphism group of the digraph is solvable and ``n-solv" otherwise;
\item {\tt Rad}:		the {\em radius}, that is, half of the length of an alternating cycle;
\item {\tt AtNo}:		the {\em intersection number}, that is,  the size of the intersection of two intersecting alternating cycles;
\item {\tt AtTy}:		the {\em attachment type}, that is: ``loose" if the attachment number is $1$, ``antipodal" if $2$, and ``tight" if  equal to the radius, otherwise ``---";
\item {\tt $|$AltCyc$|$}:	the number of alternating cycles;
\item {\tt AltExp}:	the alter-exponent;
\item {\tt AltPer}:	the alter-perimeter;
\item {\tt AltSeq}:	the alter-sequence;
\item {\tt IsGWD}:	``yes" if the digraph is generalized wreath, and ``no" otherwise.
\end{itemize}

\subsection{The data on arc-transitive $4$-GHATs}
\label{sec4.7}

The ``Census-GHAT-1k-data.csv'' file concerns arc-transitive $4$-GHATs . Each line of the file represents one of the graphs in the census, and
has nine fields, described below. Note, however, that the file does not contain the generalised wreath graphs.

\begin{itemize}
\item {\tt Name:}  the name of the graph (for example GHAT[9,1]);
\item {\tt $|$V$|$}: the order of the graph;
\item {\tt gir}:	the girth (length of a shortest cycle) of the graph;
\item {\tt bip}:	this field contains ``b" if the graph is bipartite and ``nb" otherwise;
\item {\tt CayTy}:	this field contains ``Circ" if the graph is a circulant (that is, a Cayley graph on a cyclic group), ``AbCay" if the graph is Cayley graph on an 
  abelian group, but not a circulant, and ``Cay" if it is Cayley but not on an abelian group -- it contains ``n-Cay" otherwise;
\item {\tt $|A_v|$}: the order of the vertex-stabiliser in the automorphism group of the graph;
\item {\tt $|G_v|$}:	a sequence of the orders of vertex-stabilisers of the maximal half-arc-transitive subgroups of the automorphism group -- up to conjugacy in the automorphism group;
\item {\tt solv}:	this field contains ``solve" if the automorphism group of the graph is solvable and ``n-solv" otherwise;
\item {\tt $[|$ConCyc$|]$}:  the sequence of the lengths of $A$-consistent oriented cycles of the graph (one cycle per each $A$-orbit, where $A$ is the automorphism group of the graph) -- 
		the symbols ``c'' and ``s" indicate whether the corresponding cycle is chiral or symmetric -- for example, $[4c;4c;10s]$ means there are two chiral orbits
		of $A$-consistent cycles, both containing cycles of length $4$, and one orbit of symmetric consistent cycles, containing cycles of length $10$.
\end{itemize}

\subsection{The data on $4$-HATs}

The ``Census-HAT-1k-data.csv'' file concerns $4$-HATs. Each line of the file represents one of the graphs in the census, and
has 16 fields.
The fileds {\tt $|$V$|$}, {\tt gir}, {\tt bip}, and {\tt Solv} are as in Section~\ref{sec4.7}, 
and the fields {\tt Rad, AtNo, AtTy, AltExp, AltPer} and {\tt AltSeq} are as in Section~\ref{sec4.6}.
The remaining fileds as follows:

\begin{itemize}
\item {\tt Name:}  the name of the graph (for example HAT[27,1]);
\item {\tt IsCay}:	this field contains ``Cay" if the graph is Cayley and ``n-Cay" otherwise;
\item {\tt $|G_v|$}:	the order of the vertex-stabiliser in the automorphism group of the graph;
\item {\tt CCa}:		the length of a shortest consistent cycle;
\item {\tt CCb}:		the length of a longest consistent cycle;
\item {\tt MetaCircTy}: ``$\{\}$" if the graph is not a meta-circulant; otherwise a set of types of meta-circulants that represents the graph.
\end{itemize}



\end{document}